\documentstyle[12pt,twoside, epsf]{article}

\newcommand{\ZZ}{\mbox{$Z\!\!\! Z\!$}}	

\let\\=\cr

\def\Chi{\hbox{\raise0.5ex\hbox{$\chi$}}}

\newtheorem{th}{Theorem}
\newtheorem{lem}{Lemma}
\newtheorem{cor}{Corollary}

\newtheorem{defn}{Definition}

\newtheorem{rem}{Remark}

\def\picill#1by#2(#3)
{\vbox to #2
{\hrule width #1 height 0pt depth 0pt
\vfill\epsffile{#3}}}

\let \ttorg \tt \def \tt{\ttorg \obeyspaces}

\begin{document}
\pagestyle{myheadings}

\markboth{{\sc H\"aring-Oldenburg \& Lambropoulou}}{{\sc Knot theory in handlebodies}}

\title{\Large\bf Knot theory in handlebodies}

\author{\sc Reinhard H\"aring-Oldenburg and 
Sofia Lambropoulou\thanks{A large part of this work was done as the
second author was employed at G\"ottingen Universit\"at. Both authors acknowledge with gratitude 
the research facilities offered there.} } 

\date{}

\maketitle


\begin{abstract}
We consider oriented knots and links in a handlebody of genus $g$ through appropriate braid
representatives in $S^3$, which are elements of the braid groups $B_{g,n}$. We prove a 
geometric version of the Markov theorem for braid equivalence in the handlebody, which is based on
the $L$-moves. Using this we then  prove two algebraic versions of the Markov theorem. The first one
uses the $L$-moves. The second one uses the Markov moves and conjugation in the groups $B_{g,n}$. We
show that not all conjugations correspond to isotopies.
\end{abstract}

\section{Introduction} 

A natural generalization of the classical knot theory in $S^3$ considers knots and links in more
general 3-manifolds. While topological quantum field theories provide an approach to invariants of
links in closed (i.e. compact without boundary) 3-manifolds, bounded 3-manifolds are also of
interest, since --for once-- they give rise to  closed, connected, orientable
3-manifolds. In particular,  we have on the one hand handlebodies, which give rise to 3-manifolds
via the Heegaard decomposition, and on the other hand knot complements, which give rise to
3-manifolds via the surgery technique. In \cite{LR} knots and links in knot complements and
3-manifolds are studied via braids. Here we study knots and links in handlebodies. The special case
of the solid torus is the only bounded manifold common in both categories, and its knot theory has
been studied quite extensively from various viewpoints (see \cite{T}, \cite{HP}, \cite{HK},
\cite{L1,GL,L2}, \cite{Pa}, \cite{D}, \cite{H}). Various aspects of the knot theory of a handlebody
have been studied in  \cite{P1}, \cite{P2}, \cite{S}, \cite{Li}, \cite{V}, \cite{Lie}.  Let  now
$H_g$ denote a handlebody of genus $g$. A handlebody of genus $g$ is usually defined as 
$(\mbox{\it a closed disc} \setminus \{g \, \mbox{\it open discs} \})\times I$, where $I$ is the 
unit interval. See Fig. 1.

$$\vbox{\picill7cmby1.2in(hb1)  }$$

\begin{center}
{Fig. 1. A handlebody of genus $3$.}
\end{center}
\vspace{3mm} 

\noindent Equivalently, $H_g$ can be defined as $(S^3 \setminus \mbox{\it an open tubular
neighbourhood of} \ I_g),$ where $I_g$ denotes the pointwise fixed identity braid on $g$ 
indefinitely extended strands, all meeting at the point at infinity, see Fig. 2. Thus $H_g$ {\it
may be represented in} $S^3$ by the braid  $I_g$. Now let $L$ be an oriented link in $H_g$. Then
$L$ will avoid the $g$ hollow tubes of $H_g$, and also it will not pass beyond the boundary of
$H_g$ from either end.  Equivalently,  the link $L$ in $H_g$ {\it may be  represented
unambiguously} by the mixed $(g,g)$-tangle $I_g \bigcup L$ in $S^3$, which by abuse of language we
shall call {\it mixed link} (see Fig. 2).  The subbraid $I_g$ shall be called  {\it the fixed part}
and $L$  {\it the moving part} of the mixed link. A {\it mixed link diagram} $I_g \bigcup \tilde L$
is then a diagram of $I_g \bigcup L$ projected on the plane of $I_g$, which is equipped with the
top-to-bottom direction. Note that, if we remove $I_g$ from a  mixed link we are left with an
oriented link in
$S^3$. 

$$\vbox{\picill8.4cmby2.2in(hb2)  }$$

\begin{center}
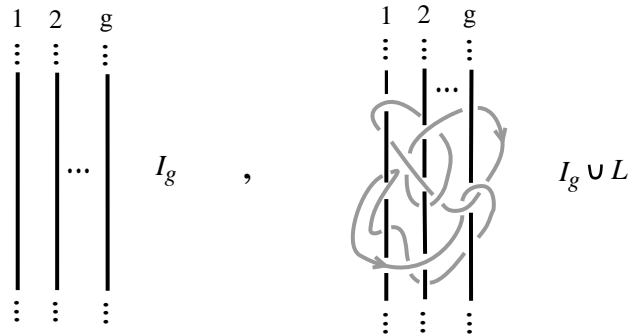
{Fig. 2. Representation of $H_g$ - a mixed link.}
\end{center}
\vspace{3mm}

In this paper we study isotopy of oriented links and equivalence of braids in $H_g$  via their 
mixed link and mixed braid representatives in $S^3$. The paper is organised as follows. In Section
2 we study knot isotopy in handlebodies combinatorially, we establish the notion of a braid in
$H_g$, and we prove that every oriented link in $H_g$ can be braided. In Section 3 we prove a
geometric version of the Markov theorem for oriented links in $H_g$ (Theorem 3) using the
so-called $L$-moves (Definition 6) and the Relative Version of Markov theorem. In Section 4 we
define the algebraic structures of braids in $S^3$ that represent braids in $H_g$ and we prove two
algebraic versions of the Markov theorem for handlebodies. The first one (Theorem 4) uses the 
$L$-moves. The second one (Theorem 5) uses a presentation of the groups $B_{g,n}$, and its
formulation resembles the classical Markov theorem for $S^3$. Only, here the Markov move (the one that
introduces a twist) has to take place anywhere on the right of the braid. Also,  as 
 we prove, not all conjugations in the groups $B_{g,n}$ induce isotopy in the
handlebody. This disproves a conjecture of A. Sossinsky, \cite{S}.  In Section 5 we discuss which
conjugations are allowed (Theorem 6).  Finally, in Section 6 we  discuss what kind of maps should be
defined on appropriate quotient algebras in order to replace the notion of a Markov trace.

\smallbreak
Parts of this paper have been presented at the AMS meeting in Buffalo, Spring '99, at the  AMS
meeting in Louisiana, Spring 2000, at the `Knots 2000' in Korea and at the KLM meeting in Siegen,
January 2001. We would like to thank the referee for several valuable comments. Also Jozef
Przytycki and Adam Sikora for very important discussions.

\section{Knots and braids in $H_g$} 

Throughout the paper the handlebody $H_g$ will be represented in $S^3$ by the braid $I_g,$ as
defined in Section 1, and a link $L$ in $H_g$ will be represented by the mixed link $I_g \bigcup L$
in $S^3$. The set-up is similar to the one  of \cite{LR},  and we will refer to
\cite{LR} for the proofs of results needed here and already established there. Otherwise, we have
tried to present our results in a self-contained manner.  All links will be assumed oriented and all
diagrams piecewise linear (PL). Whenever we say `knots'  we mean `knots and links'. Finally, we
will be thinking in terms of diagrams for both knots and braids.

\begin{defn}{\rm \ Two oriented links $L_1, L_2$ in $H_g$ are {\it isotopic} if and only
if there is an ambient isotopy of  $(S^3 \setminus I_g, \ L_1) \longrightarrow (S^3 \setminus
I_g, \ L_2)$ taking $L_1$ to $L_2$. Equivalently, $L_1$ and $L_2$ are isotopic in $H_g$ if and only
if the mixed links $I_g \bigcup L_1$ and $I_g \bigcup L_2$ are isotopic in $S^3$ by an ambient
isotopy which keeps $I_g$ pointwise fixed. } \end{defn}

\noindent In the PL category ambient isotopy is realized through a finite sequence of the
so-called {\it $\Delta$-moves} in three-space.  

\begin{defn}{\rm \ A $\Delta$-move on a link $L$ in $H_g$ is an elementary
combinatorial isotopy move (and its inverse), realized by replacing an arc of $L$ by two other arcs 
respecting orientation, and such that all three arcs span a triangle in space, the spanning surface
of which does not intersect any other arcs of $L$.  On the level of the mixed link $I_g \bigcup L$
in $S^3$, a $\Delta$-move applies only on the moving part. A 
$\Delta$-move on a mixed link diagram $I_g \bigcup \tilde L$ is the regular  projection of a
$\Delta$-move on the plane of $I_g$.  } \end{defn}

\begin{defn}{\rm \ A {\it non-critical $\Delta$-move} on a link $L$ in $H_g$ is a $\Delta$-move,
such that on its regular projection on the plane of the subbraid $I_g$ nothing critical occurs if we
remove the subbraid $I_g$. Consequently, on the level of the mixed link diagram $I_g \bigcup \tilde
L$, a non-critical  $\Delta$-move will be  a $\Delta$-move on $\tilde L$, whose spanning triangle
either does not meet any other arcs on the plane of projection (and so it is a {\it planar
$\Delta$-move} in the classical set-up) or it   meets parts of the fixed subbraid $I_g$. These last
possibilities shall be called {\it mixed isotopy  moves}, see Fig. 3. } \end{defn} 

$$\vbox{\picill13.8cmby.75in(hb3)  }$$

\begin{center}
{Fig. 3. Mixed isotopy moves.}
\end{center}
\vspace{3mm} 
  
\noindent \noindent Reidemeister \cite{Rd1}
(and Alexander, Briggs  \cite{AB})  proved that a $\Delta$-move on a link diagram in $S^3$ can break
into a finite sequence of the three local $\Delta$-moves known as `Reidemeister
moves', and of planar $\Delta$-moves, i.e. moves, whose spanning triangle  does not meet 
any other arcs on the projection plane, with their obvious  symmetries and choices of orientation.
From the above we deduce that knot isotopy in $H_g$ is realized combinatorially through the
following (compare with the relative version of Reidemeister theorem, Theorem 5.2 of \cite{LR}):  

\begin{th}[Reidemeister theorem for $H_g$]{ \ Two oriented links in $H_g$ are isotopic if and only
if any two corresponding mixed link diagrams in $S^3$ differ by a finite sequence of planar 
$\Delta$-moves, the three Reidemeister moves  and  the mixed isotopy moves (with
their obvious different choices of  orientation, crossings and direction), all of 
which apply only on the moving parts of the diagrams. 
 } \end{th}

Assuming that the strands of $I_g$ are oriented downwards, we can now define:

\begin{defn}{\rm \ A {\it geometric mixed braid} on $n$ strands, denoted $I_g \bigcup B$, is an
element of the classical braid group $B_{g+n}$, consisting of two  disjoint sets of strands, one of
which is the
 identity  braid $I_g$, whilst the other set of strands has labels `u' or `o' (for `under' or
`over')  attached  to each pair of corresponding endpoints (see Fig. 4). For the two sets of strands
we use the terms $I_g$-{\it part} for the identity  subbraid and $B_n$-{\it part} for the labelled
subbraid $B$. The reason for choosing this notation will become clear soon.  A {\it diagram of a
geometric mixed braid} is a braid diagram in the usual sense, projected on the plane of $I_g$.}  
\end{defn}

$$\vbox{\picill10.5cmby1.9in(hb4)  }$$

\begin{center}
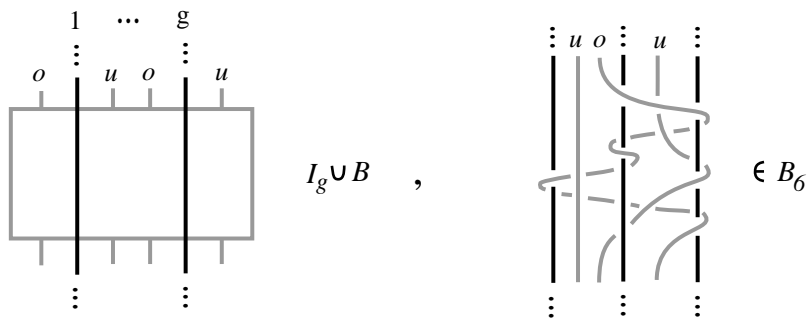
{Fig. 4. Geometric mixed braids.}
\end{center}
\vspace{3mm}

\noindent Fig. 4 illustrates an abstract geometric mixed braid enclosed in a `box', as well as an
example in $B_6$. Note that the set of geometric mixed braids on $n$ strands does not form a group,
as composition may not be well-defined.  Geometric mixed braids in $H_g$  may be visualized as
having endpoints on three different parallel planes, parallel to the plane of the paper, such that
the subbraid $I_g$ lies on the middle one, the endpoints labelled `o' lie on the front plane (the
one closest to the reader), and the endpoints labelled `u' lie on the back plane (the
 furthermost from the reader).

\bigbreak

We obtain knots from braids via the well-known closing operation adapted to our situation. So, we
have:

\begin{defn}{\rm \ The {\it closure} ${\cal C}(I_g \bigcup B)$ of a geometric mixed braid $I_g
\bigcup B$ is an operation that results in an oriented mixed link, and it is realized by joining
each pair of corresponding (slightly bent) endpoints {\it of the $B_n$-part} by a vertical segment,
either {\it over} or {\it under} the rest of the braid, according to the label attached to these
endpoints (see Fig. 5 for an example).} \end{defn}

$$\vbox{\picill4cmby1.9in(hb5)  }$$

\begin{center}
{Fig. 5. Closure of a geometric mixed braid.}
\end{center}
\vspace{3mm} 
 
\noindent Note that the strands of $I_g$ do not participate in the closure operation, that's why they
are assumed to be infinitely extensible. Besides, the labelling `u' or `o' for corresponding 
endpoints in Definition 4 is precisely an instruction on how to perform the closure. Different
choices of labels will yield in general non-isotopic links in $H_g$, as the example in Fig. 6 
illustrates. We return to this example in the discussion before Fig. 19.  

$$\vbox{\picill13.4cmby1.45in(hb6)  }$$

\begin{center}
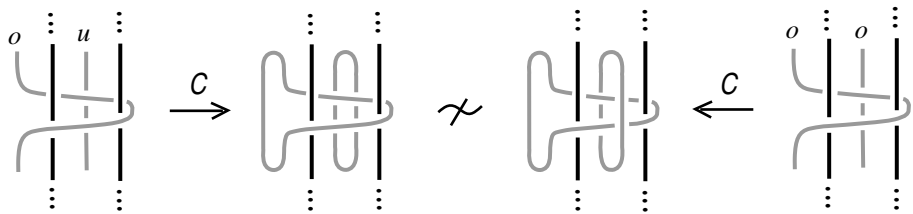
{Fig. 6. Different labels yield non-isotopic links.}
\end{center}
\vspace{3mm} 
 
\begin{rem}{\rm \  Let $M$ denote  the complement of the
$g$-unlink or a connected sum of $g$ lens spaces of  type $L(p,1)$. Then braids in $M$ can be also
represented in $S^3$ by unlabelled geometric mixed braids with $I_g$ as a fixed subbraid (cf.
\cite{LR}). Note that in both  $H_g$ and $M$, if we remove $I_g$ from a mixed braid, we are left
with a  braid in $S^3$. This will be a labelled braid in the case of $H_g$. But this is equivalent to
the familiar unlabelled picture of a classical braid, since a closing arc labelled `o' can slide
freely over to the side and then to the back of the braid, thus aquiring the label `u' (see Fig. 7).
This isotopy is the reason that mixed braids in $M$ are not labelled, since  in the set-up of
\cite{LR} $I_g$ participates also in the closure of the braid (contrary to  $H_g$).  
 } \end{rem}

$$\vbox{\picill12.2cmby1.45in(hb7)  }$$

\begin{center}
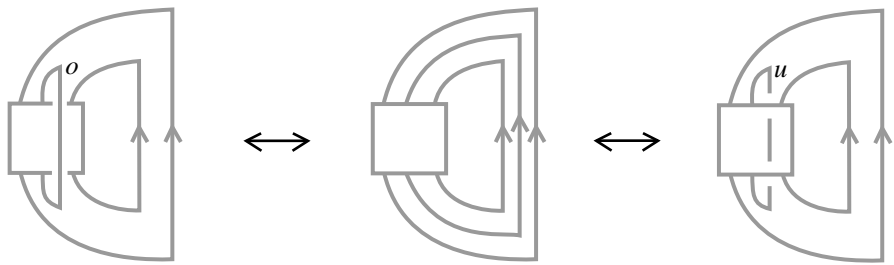
{Fig. 7. The `under - over' interchange.}
\end{center}
\vspace{3mm} 

Conversely to the closure of braids, mixed links may be braided, so that if we start with a 
mixed link, do braiding and then take closure, we obtain a mixed link isotopic to the original one.
Indeed, we have: 

\begin{th}[Alexander theorem for $H_g$]{ \ An oriented mixed link $I_g \bigcup L$ in $H_g$ may be
 braided to a geometric mixed braid, the closure of which is isotopic to $I_g \bigcup L$. 
 } \end{th}

\bigbreak
\noindent {\em Proof.} \  We apply the braiding algorithm of \cite{LR} on a  diagram
$I_g \bigcup \tilde L$ of the PL mixed link $I_g \bigcup L$.  
By general position $I_g \bigcup \tilde L$ contains no horizontal arcs with
respect to the height function.  The idea of the braiding is on the one hand to keep the arcs of the
diagram that are oriented downwards with respect to the height function, and on the other hand to
eliminate the ones that go upwards and produce instead  braid strands.  We call these arcs {\it 
opposite arcs}.  Now, the point is that the subbraid  $I_g$ will not be touched by the algorithm,
so the opposite arcs will be arcs of the link $L$.  The elimination of the opposite arcs is based on
the following: If we run along an opposite arc we are likely to meet a succession of overcrossings
and undercrossings. We subdivide (marking with points) every opposite arc into smaller -- if
necessary -- pieces, each containing crossings of only  {\it one}  type;  i.e. we may have:

$$\vbox{\picill4.5inby.9in(hb8)  }$$

\begin{center}
{Fig. 8. }
\end{center}
\vspace{3mm} 

\noindent We call the resulting pieces {\it  up-arcs}, and we label every 
 up-arc with an  {\it `o'} resp. {\it `u'} according as it is the  {\it over} resp. {\it under}  arc
of a crossing (or some crossings). If it is a  {\it free up-arc}
(one that contains no crossings), then we have a choice whether  to label it {\it `o'} 
or  {\it `u'}. The idea is to eliminate the  opposite arc by eliminating its up-arcs one
by one and create instead a pair of braid strands for each up-arc. 

  Let now $P$ be the top vertex of the up-arc $QP$ (see Fig. 9). Associated to $QP$ is the {\it
sliding   triangle} $T(P)$, which is a special case of a triangle
needed for a $\Delta$--move; it is right-angled with hypotenuse  $QP$
and with the right angle lying below the up-arc. Note that,  if $QP$ is itself vertical, then
$T(P)$  degenerates into the arc $QP$.  We say that a sliding triangle is of type {\it over} or
{\it under} according to the label of the up-arc it is associated with. (This implies that there
may be triangles of the same type lying one on top of the other.)

 The germ of our braiding process is this. Suppose for definiteness that $QP$  is of
type  {\it over}. Then we cut $QP$ at $P$  and we pull the two ends, the top upwards and the lower
downwards, and both {\it over} the rest of the diagram, so as to create a pair of corresponding
strands of the anticipated braid  (see Fig. 9).  Finally, we perform  a $\Delta$-move across the
sliding triangle $T(P)$. By general position the resulting diagram will be regular and $QQ'$ may be
assumed to slope slightly downwards. If $QP$ were {\it under} then the pulling of the two ends would
be  {\it under} the rest of the diagram.
Note that the effect of these two operations has been  to replace the up-arc $QP$ by three arcs none
of which is an up-arc, and the two of them being corresponding braid strands. Therefore we now have
fewer up-arcs. 

$$\vbox{\picill4.5inby1.4in(hb9)  }$$ 

\begin{center}
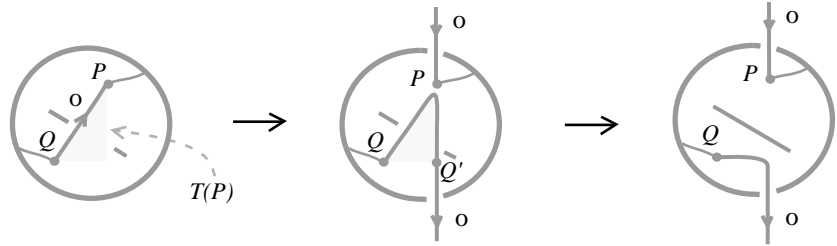
{Fig. 9. The germ of the braiding.}
\end{center}
\vspace{3mm} 

  For each up-arc that we eliminate, we label the corresponding
end strands `o' or `u' according to the label of their up-arc. (As already noted, in \cite{LR} this
labelling was not needed.)  After eliminating all up-arcs we  obtain a geometric mixed braid,
denoted  ${\cal B}(I_g \bigcup L),$ the closure of which is obviously isotopic to $I_g \bigcup L$. 
Indeed,  from Definition 5, the closing arc of  two corresponding endpoints of the braid is
precisely a stretched version of the initial up-arc, since it bears the same label. 
 $\hfill \Box $  

\noindent The proof of Theorem 2 is analogous to the one in Section 3 of \cite{LR}. We have repeated it
here for the sake of completion.

\section{Geometric Markov theorem for $H_g$}

As in classical knot theory, the next consideration is how to characterize geometric mixed braids
that induce via closure isotopic links in $H_g$. For this purpose we need to recall the $L$-moves
between braids. These were introduced in \cite{LR}, and they generalize geometrically the Markov
moves in the following sense. An {\it $L_o$-move} resp. {\it $L_u$-move} on a braid consists of
cutting an arc open and splicing into the broken strand new arcs to the top and bottom, both {\it
over} resp. {\it under} the rest of the braid (see Fig. 10 for the case of $H_g$). As remarked
in \cite{LR}, using a small braid isotopy, a braid $L$-move can be equivalently seen with a crossing
(positive or negative) formed (see Fig. 11 for $H_g$). Therefore, a geometric Markov move in a braid,
that introduces a crossing in the bottom right  position, is a special case of an $L$-move.
$L$-moves and braid isotopy generate an equivalence relation on braids called {\it $L$-equivalence}.
It was shown in \cite{LR} that $L$-equivalent classes of braids are in bijective correspondence with
isotopy classes of oriented links in $S^3$,  the bijection being  induced by `closing' the braid.
Modified slightly, $L$-moves in a handlebody are defined as follows.

\begin{defn}[Geometric $L$-moves in $H_g$ ]{\rm \ Let $I_g \bigcup B$ be a braid in $H_g$ and $P$ a
point of an arc of the subbraid $B$, such that $P$ is not vertically aligned with any  crossing.
Doing a {\it geometric $L$-move}  at $P$ means to perform the following
operation:  cut the arc at $P$, bend the two resulting smaller arcs slightly apart  by a small
isotopy and stretch them vertically, the upper downwards and the lower upwards, and both {\it over} 
or {\it under} all other arcs of the diagram, so as to introduce two new corresponding strands with
endpoints on the vertical line of $P$,  labelled `o' or `u' according to the stretching. Stretching
the new strands over  will give rise to a {\it geometric $L_o$-move} and under to an {\it
geometric $L_u$-move}. Undoing an $L$-move is defined to be the reverse operation. Also in this
set-up, two geometric mixed braids in $H_g$ that differ  by an $L$-move shall be called {\it
$L$-equivalent}.  }\end{defn}

$$\vbox{\picill13.5cmby1.85in(hb10)  }$$

\begin{center}
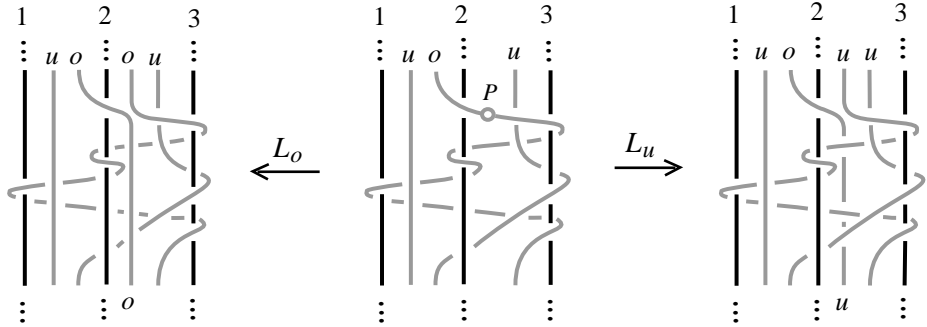
{Fig. 10. The two types of $L$-moves in $H_g$.}
\end{center}
\vspace{3mm} 

\noindent Fig. 10 illustrates an example of a geometric $L_o$-move and a geometric $L_u$-move at
the same point of a geometric mixed braid, whilst Fig. 11 illustrates an abstract geometric
$L_o$-move and the  crossing it introduces in the braid box. 

$$\vbox{\picill8inby1.65in(hb11)  }$$

\begin{center}
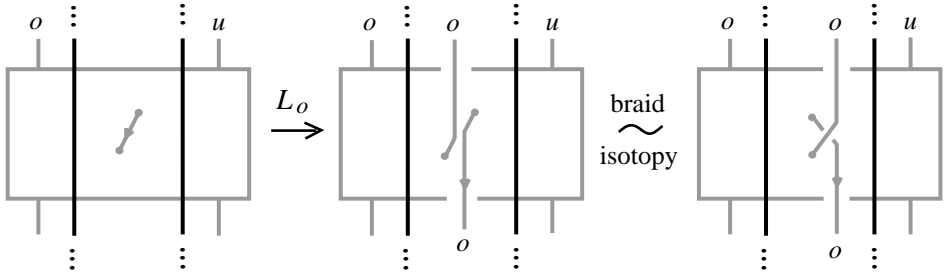
{Fig. 11. An $L$-move introduces a crossing.}
\end{center}
\vspace{3mm} 
 
\begin{rem}{\rm \ $L$-equivalent geometric mixed braids have isotopic
closures, since the labels we give to the new endpoints after performing an $L$-move on a
mixed braid agree with the type of the $L$-move. So closure is compatible with the
$L$-move, and it  corresponds to introducing a twist in the mixed link.   
 } \end{rem}

We are now in a position to state the following. 

\begin{th}[Geometric version of Markov theorem for $H_g$]{ \ Two oriented links in $H_g$ are 
isotopic iff any two corresponding geometric mixed braids differ by a finite sequence of $L$-moves
and isotopies of geometric mixed braids.
 } \end{th}

\noindent {\em Proof.} \  Let ${\cal B}$ denote the braiding map of Theorem 2, let $I_g \bigcup
\tilde L$ be a mixed link diagram in $H_g$, and let $I_g \bigcup B = {\cal B}(I_g \bigcup \tilde
L)$. By Theorem 2, \ ${\cal C}\circ {\cal B}(I_g \bigcup \tilde L)$ is isotopic to 
$I_g \bigcup \tilde L$. Further, ${\cal B}\circ {\cal C}(I_g \bigcup B) = I_g \bigcup
B$. This follows from Definition 5 and from the fact that if we braid the closing arcs of a mixed
braid, $I_g \bigcup B$ say, each closing arc will give rise to one pair of corresponding strands,
so we obtain again the braid $I_g \bigcup B$.

\bigbreak
\noindent We now consider the liftings ${\tilde {\cal B}}$ and ${\tilde {\cal C}}$ of the maps
${\cal B}$ and ${\cal C}$ on isotopy classes of link diagrams and on $L$-equivalent classes of
geometric mixed  braids respectively. We will show that ${\tilde {\cal C}}$ is a bijection with
inverse  ${\tilde {\cal B}}$. It follows from Remark 2 that ${\tilde {\cal C}}$ is well-defined.
Thus, from the observations above, it only remains to show that ${\tilde {\cal B}}$ is also
well-defined, that is to show that geometric mixed braids corresponding to isotopic mixed links are
$L$-equivalent.  For this we apply the: 

\bigbreak

\noindent {\bf Relative Version of Markov theorem} (Theorem 4.7 of \cite{LR}) \ {\it Let $L_1$,
$L_2$ be oriented link diagrams in $S^3$, both containing a common braided portion $B$. Suppose that
there is an isotopy of $L_1$ to $L_2$ which finishes with a homeomorphism fixed on $B$. Suppose
further that $B_1$ and $B_2$ are braids obtained from our braiding process applied to $L_1$ and
$L_2$ respectively, both containing the  common braided portion $B$. Then $B_1$ and $B_2$ are
$L$--equivalent by moves that do not affect the braid $B$. }

\bigbreak

\noindent Here $I_g$ plays the role of the the common subbraid $B$, which,  by
Definition 3 and by Theorem 1,  remains fixed throughout an isotopy of two mixed link diagrams $I_g
\bigcup \tilde L_1$ and $I_g \bigcup \tilde L_2$. Further, the braiding ${\cal B}$ keeps $I_g$ fixed
in the corresponding geometric mixed braids, $I_g \bigcup B_1$ and $I_g \bigcup B_2$, say. Thus, the
relative version of Markov theorem guarantees that $I_g \bigcup B_1$ and $I_g \bigcup B_2$ are
$L$-equivalent by  $L$-moves that do not affect $I_g$. But this is precisely the definition of 
$L$-moves in $H_g$ (recall Definition 6). The only difference from $S^3$ is that 
here we attach labels to the corresponding strands of each $L$-move according to its type. In
$S^3$ this was not needed.  $\hfill \Box  $

\section{Algebraic versions of Markov theorem }

In order to  construct invariants of knots in the handlebody using the braid approach  we must 
translate  Theorem 3 into algebra (see for
example \cite{J} for $S^3$ and \cite{L2} for the solid torus). For this we need first to introduce the
braid groups $B_{g,n}$. 

\begin{defn}{\rm \ An {\it algebraic mixed braid} on $n$ strands is
an element of the braid group $B_{g+n}$ consisting of two  disjoint sets of strands, such that
{\it the first} $g$ {\it strands} constitute the identity braid $I_g$. We denote algebraic mixed
braids in the same way as the geometric mixed braids. 
  } \end{defn}

\noindent Fig. 12 suggests two ways for depicting abstractly algebraic mixed braids, and it gives a
concrete example of an algebraic mixed braid on three strands.

$$\vbox{\picill13cmby1.95in(hb12)  }$$

\begin{center}
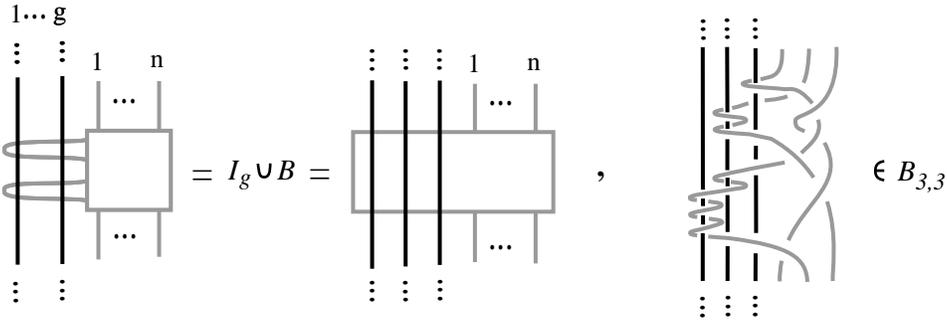
{Fig. 12. Algebraic mixed braids.}
\end{center}
\vspace{3mm}

\noindent We shall see that an algebraic mixed braid is a special case of a geometric mixed braid.
Clearly, it is a special case of an {\it unlabelled} geometric mixed braid. We say that an algebraic
mixed braid {\it is made geometric} if we attach arbitrary labels `u' or `o' at its corresponding
endpoints. Note that this is an ambiguous process. 

\begin{defn}{\rm \ The {\it closure} of an algebraic mixed braid $I_g \bigcup B,$ denoted $I_g
\bigcup \widehat{B},$ is defined by joining each pair of the (slightly bent) corresponding endpoints
of the $B_n$-part by a vertical segment (see left illustration of Fig. 13). 
  } \end{defn} 

$$\vbox{\picill10.5cmby1.95in(hb13)  }$$

\begin{center}
{Fig. 13. Closure of an algebraic mixed braid.}
\end{center}

\begin{rem}{\rm \ If we consider an algebraic mixed braid $I_g \bigcup B$ made
geometric, its closure  ${\cal C}(I_g \bigcup B)$ is isotopic to $I_g \bigcup \widehat{B}$, no
matter what labels we used for the $B_n$-part, since the closing arcs can be stretched and can slide
freely over to the right-hand side of the braid (see  Fig. 13). {\it This shows that algebraic mixed
braids are indeed special cases of geometric mixed braids, for which labels are superfluous. }
 } \end{rem}

\noindent Conversely, geometric mixed braids can be made algebraic. Indeed, the operation `closure' 
is an equivalence relation in the set of geometric mixed braids, and we have: 

\begin{lem}{ \ Every geometric mixed braid may be represented by an algebraic mixed braid with
isotopic closure. 
 } \end{lem}

\noindent {\em Proof.} \  Pull each pair of corresponding endpoints of the geometric mixed braid 
$I_g \bigcup B$ to the right side of $I_g$ {\it over} or {\it under} the strands of $I_g$ according
to its label, starting from the rightmost pair, and respecting the position of the endpoints.
Schematically: 

$$\vbox{\picill4.7inby3.35in(hb14)  }$$

\begin{center}
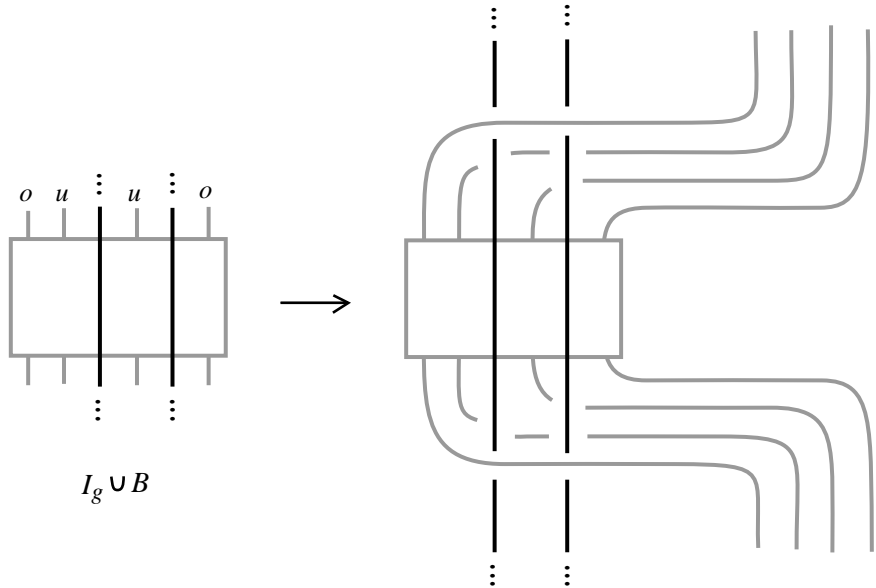
{Fig. 14. Algebraization of a geometric mixed braid.}
\end{center}
\vspace{3mm} 

\noindent We thus obtain unambiguously an algebraic mixed braid. We denote this last step of the
braiding algorithm by ${\cal A}$, and we say that through ${\cal A}$ a geometric mixed braid is
made algebraic. Now, ${\cal C}(I_g \bigcup B)$ is isotopic to the closure of the algebraic mixed
braid ${{\cal A}(I_g \bigcup B)}$. To see this we choose as labels of the algebraic mixed braid
${{\cal A}(I_g \bigcup B)}$ the initial labels of the geometric mixed braid $I_g \bigcup B$. Then
the closures of the two geometric mixed braids are isotopic, and, by Remark 3 above, the assertion
is proved.       $\hfill \Box  $

As an example, the algebraization of the two geometric braids of Fig. 6 are illustrated in Fig. 21.

\bigbreak

  The sets of algebraic mixed braids on $n$ strands, denoted $B_{g,n}$, have been treated in
\cite{L3}. It is shown there that these are the underlying braid structures for studying knots
in a handlebody, in knot complements and in closed, connected, orientable 3-manifolds. Moreover,
they form subgroups of the groups $B_{g+n}$ with operation the usual concatenation, and with
presentation: 

\[ B_{g,n} = \left< \begin{array}{ll}  \begin{array}{l} 
a_1, \ldots, a_g,  \\ 
\sigma_1, \ldots ,\sigma_{n-1}  \\
\end{array} &
\left|
\begin{array}{l} \sigma_k \sigma_j=\sigma_j \sigma_k, \ \ |k-j|>1   \\ 
\sigma_k \sigma_{k+1} \sigma_k = \sigma_{k+1} \sigma_k \sigma_{k+1}, \ \  1 \leq k \leq n-1  \\
{a_i} \sigma_k = \sigma_k {a_i}, \ \ k \geq 2, \   1 \leq i \leq g,    \\ 
 {a_i} \sigma_1 {a_i} \sigma_1 = \sigma_1 {a_i} \sigma_1 {a_i}, \ \ 1 \leq i \leq g  \\
 {a_i} (\sigma_1 {a_r} {\sigma_1}^{-1}) =  (\sigma_1 {a_r} {\sigma_1}^{-1})  {a_i}, \ \ r < i.  
\end{array} \right.  \end{array} \right>,  \]
 
\noindent where the generators of $B_{g,n}$ may be represented geometrically by: 

$$\vbox{\picill12.2cmby1.5in(hb15)  }$$

\begin{center}
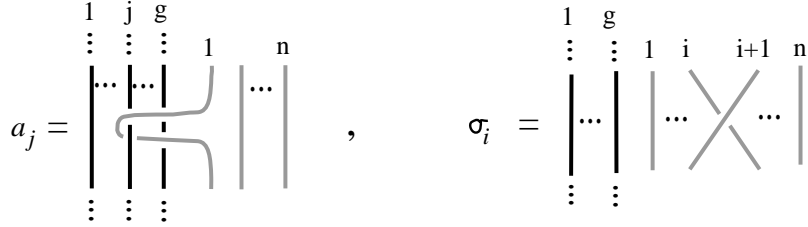
{Fig. 15. The generators of $B_{g,n}$.}
\end{center}
\vspace{3mm}

\noindent Let $B_{g,\infty} := \bigcup_{n=1}^{\infty} B_{g,n}$ denote the disjoint union of all
braid groups $B_{g,n}$ (not the inductive limit). Proceeding towards the algebraization of
Theorem 3 we define:

\begin{defn}{\rm \ An {\it algebraic $L$-move}  is a geometric $L$-move
between elements of  $\bigcup_{n=1}^{\infty} B_{g,n}$, i.e. an $L$-move that preserves the
group structure of the algebraic mixed braids.} \end{defn}

\noindent  It follows from Remark 3 that algebraic $L$-moves do not need  the labels `o'
and `u'. In some illustrations we keep the labels for the sake of clarity. An algebraic $L$-move in
a braid $\alpha \in B_{g,n}$ has the following algebraic expressions, depending on its type. These
are easily derived, as Fig. 16 shows.  

\bigbreak

$L_o$-{\it type:} \  $\alpha=\alpha_1\alpha_2 \ \sim \ 
\sigma_i^{-1}\ldots \sigma_n^{-1} \alpha_1\sigma_{i-1}^{-1}\ldots
\sigma_{n-1}^{-1}\sigma_n^{\pm 1}\sigma_{n-1} \ldots \sigma_i
\alpha_2\sigma_n \ldots \sigma_i,$

\bigbreak

$L_u$-{\it type:} \  $\alpha=\alpha_1\alpha_2 \ \sim \ 
\sigma_i\ldots \sigma_n \alpha_1\sigma_{i-1}\ldots
\sigma_{n-1}\sigma_n^{\pm 1}\sigma_{n-1}^{-1}\ldots\sigma_i^{-1}
\alpha_2\sigma_n^{-1}\ldots\sigma_i^{-1}.$

$$\vbox{\picill15.5cmby2in(hb16)  }$$

\begin{center}
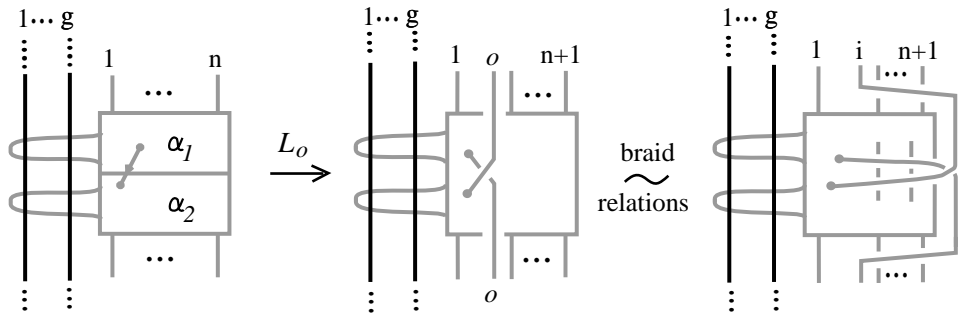
{Fig. 16. Algebraic expression of an algebraic $L_o$-move.}
\end{center}
\vspace{3mm} 
 
\begin{lem}{ \ Consider a geometric mixed braid containing a geometric $L$-move, which is made
algebraic. Then the $L$-move is turned into an algebraic $L$-move. 
  } \end{lem}

\noindent {\em Proof.} \  Since the type of a geometric
$L$-move agrees with the label of its endpoints, by pulling the endpoints to the right  the
crossing of the  $L$-move slides over by a braid isotopy. Schematically:  

$$\vbox{\picill13.5cmby2.6in(hb17)  }$$

\begin{center}
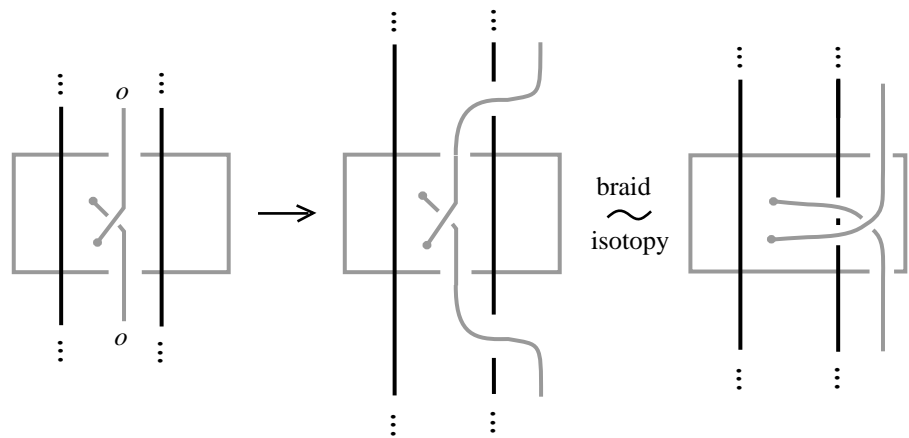
{Fig. 17. Sliding a geometric $L_o$-move to the right.}
\end{center}
\vspace{3mm} 
 


\noindent The case of a geometric $L_u$-move is completely analogous. Here the pulling takes place
under the braid, so the crossing of the geometric $L_u$-move slides along to the right to form an
algebraic $L_u$-move. 
 $\hfill \Box $

\bigbreak
Now we can state the following:
 
\begin{th}[First algebraic version of Markov theorem for $H_g$]{ \ Two oriented links in
$H_g$ are isotopic iff any two corresponding algebraic mixed braids differ by a finite sequence of 
 algebraic $L$-moves and the braid relations in $\bigcup_{n=1}^{\infty} B_{g,n}$.
 } \end{th}

\noindent {\em Proof.} \ It follows from Theorem 3 and Lemma 2.
 $\hfill \Box  $

\begin{rem}{\rm \ Theorem 4 is a rephrasing of Theorem 3 in an algebraic set-up. One could omit
Theorem 3 and prove Theorem 4 directly using the Relative Version of Markov theorem, after
incorporating in the braiding  algorithm  ${\cal B}$ the last algebraization step ${\cal A}$. We
decided to separate the two results, so as to stress the passage from the geometric to the algebraic
set-up and the results that are valid in each one.
 } \end{rem}

\noindent  In order to look for Markov functionals on $B_{g,\infty}$, so as to 
construct link invariants in $H_g,$ we further prove:  

\begin{th}[Second algebraic version of Markov theorem for $H_g$]{ \ \ Two oriented links in  $H_g$ are
isotopic iff any two corresponding algebraic mixed braids  differ by a finite sequence of the
following moves: 
\begin{enumerate}
\item Markov move: \ \ $\beta_1\beta_2\sim\beta_1\sigma_n^{\pm 1}\beta_2,
   \quad\beta_1,\beta_2\in B_{g,n}$
\item Markov conjugation: \ \ ${\sigma_i}^{-1}\beta \sigma_i \sim \beta,
   \quad\beta,\sigma_i\in B_{g,n}$
\end{enumerate}
 } \end{th}

\noindent {\em Proof.} \ The two types of moves are illustrated in Figs. 18a and 18b. It is
easy to see that both do not leave the isotopy class of the link. In fact, the first one is
simply a special case of an algebraic $L$-move that takes place at the rightmost part of the
algebraic mixed braid, whilst the second one clearly induces isotopy via closure, as defined
in Definition 8. The other direction is shown by reducing to Theorem 4. Indeed, an
algebraic  $L$-move can be realized by a finite sequence of the above moves, as it follows
clearly from the algebraic expressions of the two types of algebraic $L$-moves (recall Fig.
16). $\hfill \Box  $

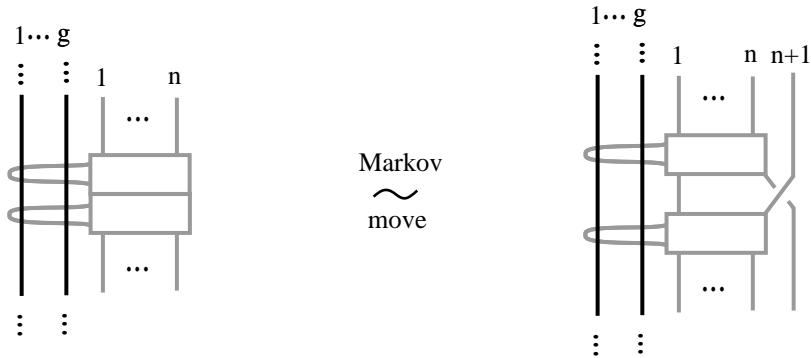
$$\vbox{\picill12cmby2in(hb18a)  }$$

\begin{center}
{Fig. 18a. The Markov move in $H_g$.}
\end{center}
\vspace{3mm} 

$$\vbox{\picill12cmby1.7in(hb18b)  }$$

\begin{center}
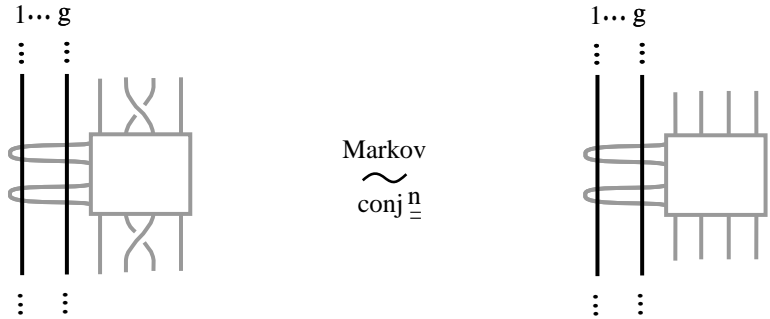
{Fig. 18b. Markov conjugation in $H_g$.}
\end{center}
\vspace{3mm} 

\noindent Algebraic mixed braids that are equivalent in the context of Theorem 4 or Theorem 5 shall
be called {\it Markov equivalent}.  A remark is now in order. 

\begin{rem}{\rm \ In the classical case in $S^3$ the braid move  
$\alpha_1\sigma_n^{\pm 1}\alpha_2 \sim \alpha_1 \alpha_2$ is equivalent to the move 

 $$\alpha\sigma_n^{\pm 1}\sim \alpha,$$
 where $\alpha = \alpha_1 \alpha_2 \in B_n$. This is still true in the case of a solid torus
(see \cite{L1}, \cite{H}). To see this think that the infinite strand of a solid torus may be closed
at the point at infinity, so any loop can conjugate with no obstruction. In \cite{H},  Lemma 39, it
is shown how to commute a loop from the bottom to the top of the braid without closing the infinite
strand. But in a handlebody of genus greater than one this is not the case any more. Here the braid
word  $\alpha_2$ may contain more than one of the $g$ generators $a_i$ of the braid group  (recall
Fig. 15). This is discussed in detail in the next section.
 } \end{rem}

\section{On hidden conjugations}

There are two natural questions arising now: 
\begin{itemize}
\item[(1)]  are there any `hidden' conjugations involving the
generators $a_i$, which preserve the isotopy class of the closure of a mixed
braid (even though the strands of $I_g$ do not participate in the closure)? 

\item[(2)] if yes, are all conjugations `allowed'? 
\end{itemize}
\noindent Before answering we need to introduce another notion. 

\begin{defn}{\rm \ A {\it loop} in $B_{g,n}$ is a word of the  form $b_i := a_i a_{i+1}
\cdots a_g$ or its inverse, for $i < g$, and a {\it maximal loop} the word $b_1 := a_1 a_2 \cdots
a_g$ or its inverse (see Fig. 19 for illustrations). A maximal loop shall be denoted by $\omega$. }
\end{defn} 

$$\vbox{\picill13cmby1.35in(hb19)  }$$

\begin{center}
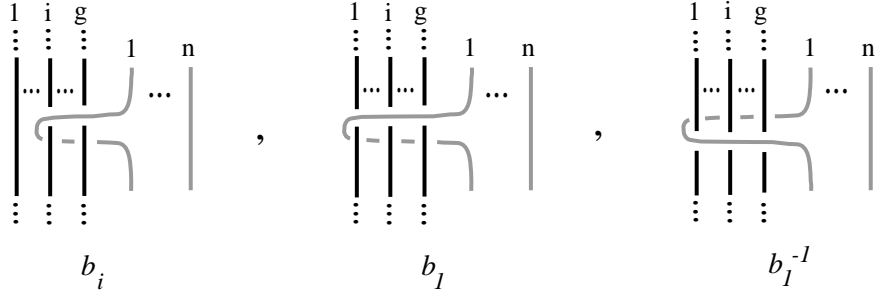
{Fig. 19. A loop and the two maximal loops.}
\end{center}
\vspace{3mm}

\noindent The answer to the first question is  positive.
Indeed, we have: 

\begin{lem}{ Let  $\alpha \in B_{g,n}$ be arbitrary and let $\omega$ be a maximal loop. Then the
braids $\alpha \, \omega$ and  $\omega \, \alpha$ are Markov equivalent.
 } \end{lem}

\noindent {\em Proof.} Fig. 20 demonstrates that, using algebraic $L$-moves and conjugation by a 
$\sigma_1$, the given algebraic mixed braids are Markov equivalent by Theorem 5. Thus their
closures are isotopic.  $\hfill \Box  $ 

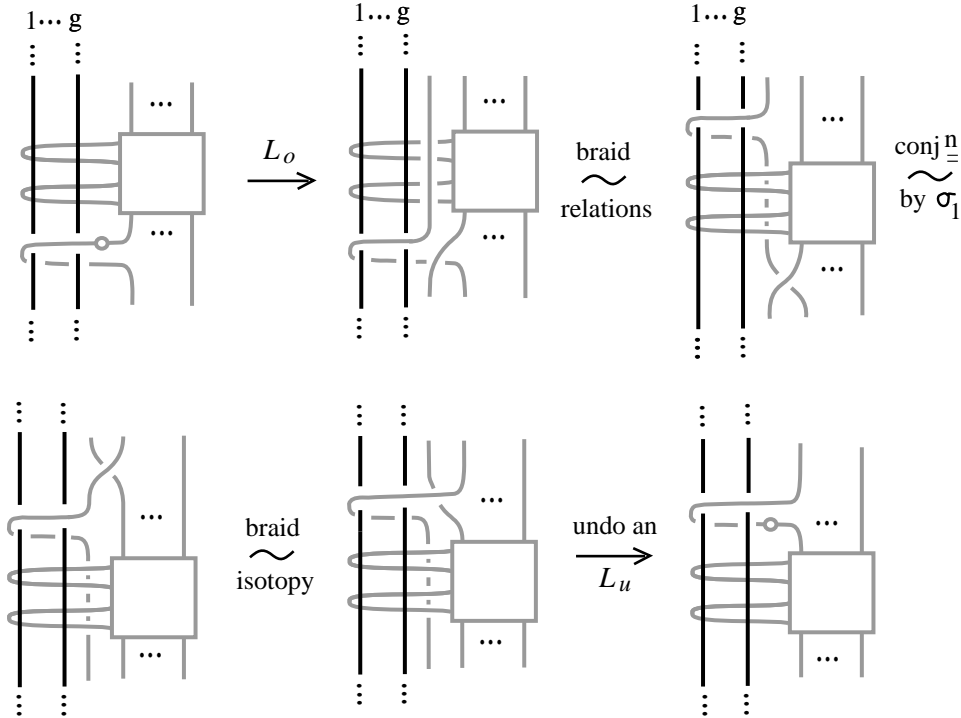
$$\vbox{\picill15cmby1.5in(hb20a)  }$$

$$\vbox{\picill15cmby1.7in(hb20b)  }$$

\begin{center}
{Fig. 20. The proof of Lemma 3.}
\end{center}
\vspace{3mm}

\begin{rem}{ \rm  An alternative proof of Lemma 3 would be to take the closures of ${\omega}^{-1} \,
\alpha \, \omega$ and $\alpha$ and to observe that a closing arc of ${\omega}^{-1} \, \alpha \,
\omega$ can be dragged around to the left and all the way round to the position $\alpha.$
 } \end{rem}

\noindent The answer to the second question is negative, as the example below shows: The
 two algebraic mixed braids of Fig. 21 are the algebraizations of the geometric mixed braids of
Fig. 6, thus their closures are the two non-isotopic mixed links of Fig. 6. 

$$\vbox{\picill9.5cmby1.8in(hb21)  }$$

\begin{center}
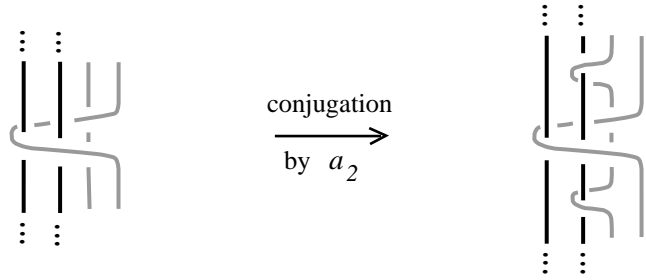
{Fig. 21. Two non-isotopic conjugate algebraic mixed braids.}
\end{center}
\vspace{3mm}

\noindent  Now a third question arises: 
\begin{itemize}
\item[(3)]   Can we list all the conjugations that are allowed with respect to isotopy?  
\end{itemize}
\noindent In order to answer this question we give first another 
presentation of the braid group $B_{g,n}$ with the $b_i$'s as generators. This presentation is
easily derived from the one with the generators $a_i$ given in the previous section using that 
  $a_i = b_i {b_{i+1}}^{-1}.$  

\[ B_{g,n} = \left< \begin{array}{ll}  \begin{array}{l} 
b_1, \ldots, b_g,  \\ 
\sigma_1, \ldots ,\sigma_{n-1}  \\
\end{array} &
\left|
\begin{array}{l} \sigma_k \sigma_j=\sigma_j \sigma_k, \quad |k-j|>1   \\ 
\sigma_k \sigma_{k+1} \sigma_k = \sigma_{k+1} \sigma_k \sigma_{k+1}, \quad  1 \leq k \leq n-1  \\
{b_i} \sigma_k = \sigma_k {b_i}, \quad k \geq 2, \   1 \leq i \leq g    \\ 
 {b_i} \sigma_1 {b_r} \sigma_1 = \sigma_1 {b_r} \sigma_1 {b_i}, \quad  r \leq i  
\end{array} \right.  \end{array} \right>.  \]
 
\noindent  It is important to understand that {\it conjugation by some}
$b_i$ of an algebraic mixed braid {\it is equivalent to changing the labels} of some pair of
corresponding endpoints in a related geometric mixed braid. If, for example, the two corresponding
endpoints of the geometric mixed braid lie to the left of all strands of $I_g$, then by a braid isotopy
part of which is absorbed inside the braid box, the geometric mixed braid can look like one of the
middle pictures of Fig. 22. Thus, change of labels corresponds to conjugating the algebraization of the
`o' braid (resp. the `u' braid)  by the loop
$b_1$ (resp.
${b_1}^{-1}$), as Fig. 22 demonstrates. 

$$\vbox{\picill15cmby2.5in(hb22)  }$$

\begin{center}
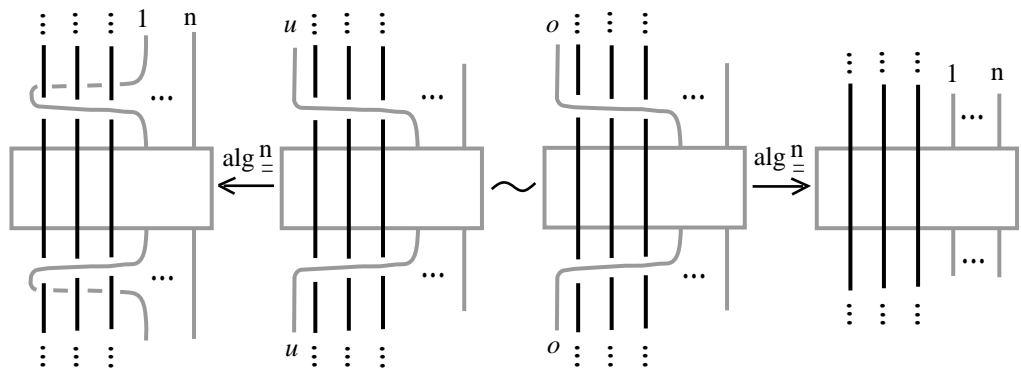
{Fig. 22. Conjugation by $b_1$ corresponds to allowed change of labels.}
\end{center}
\vspace{3mm} 

\noindent But these conjugations are allowed as we proved in Lemma 3. (To see the isotopy on the
level of the geometric mixed braids look at the closing arc of the left middle braid of Fig.
22. This can pass from the `u' position to the `o' position without any obstruction from
the braid.) Thus, change of labels in this case reflects isotopy between the
closures of the two geometric mixed braids.  

  Let, now, the two corresponding endpoints of a geometric mixed braid lie between the $i$th and
$(i+1)$st strand of  $I_g,$ for  $i\neq 1.$ Then, with similar reasoning as above, change of labels
corresponds to conjugating the `o' braid (resp. the `u' braid) by the loop  $b_{i+1}$ (resp.
${b_{i+1}}^{-1}$), see Fig. 23. And conversely, two algebraic mixed braids that are conjugate by a $b_i,$
with $i\neq 1,$ can be seen as the algebraizations of two geometric mixed braids which differ only by
the labels of one pair of corresponding endpoints. But such a change of labels does not reflect
isotopy. To see
this, think of the infinitely extended strands of
$I_g$ joining at the point at infinity. Then, the closing arc of the geometric `u' braid would have to
cross the point at infinity in order to come to the `o' position.  Consequently  conjugations by the 
$b_i$'s for  $i\neq 1$ are not allowed,  except for some obvious special cases of disconnected diagrams,
which would then imply that the knot can also live in a handlebody of smaller genus. 

$$\vbox{\picill15cmby2.3in(hb23)  }$$

\begin{center}
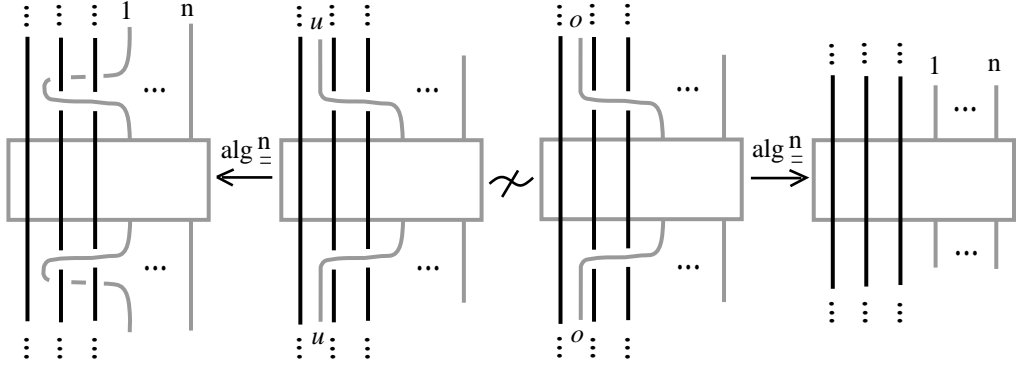
{Fig. 23. Conjugation by $b_2$ corresponds to non-allowed change of labels.}
\end{center}
\vspace{3mm}

Finally, the answer to the third question lies in the following result.

\begin{th}{ \ Conjugation by a maximal loop $\omega$ is the only conjugation by words
in the $b_i$'s, which preserves the isotopy class in $H_g$ of the closure of {\em any} braid in
$B_{g,n}.$ That is, for any word $\beta \in B_{g,n}$ in the  $b_i$'s, different from ${b_1}^n,$ for
$n \in \ZZ,$ there is an  $\alpha \in B_{g,n},$ such that $\alpha \, \beta$ and $\beta \, \alpha$
have non-isotopic closures. 
} \end{th} 

\noindent {\em Proof.}  Counter-examples of the required kind exist in the one-strand braid group
$B_{g,1}$. Indeed, assume the theorem were false, i.e. the closures of $\beta \, \alpha \,
\beta^{-1} \in B_{g,1}$ and $\alpha \in B_{g,1}$ are isotopic, for all $\alpha, \, \beta \in
B_{g,1}$ with $\alpha$ a word in the $b_i$'s and $\beta$ some $b_r.$ From Theorem 4 we know that
these two braids are related by braid isotopies and algebraic $L$-moves. Since $B_{g,1}$ is the group
generated by the $b_i$'s, and this is a free group, $\beta$ cannot be commuted through $\alpha$.
Hence we have to invoke the $L$-moves. These introduce some $\sigma_i$'s  but they do not change the
order of the $b_i$'s. According to the relations in the braid group, this can only be done if the
condition $r \leq i$ is satisfied. But this is true always, only if $r=1$. Therefore, $\beta$ has to
be $b_1$ or its inverse. Thus, conjugation by $b_r$ for $r \geq 2$ cannot be realized in the generic
case. Topologically, this corresponds precisely to crossing the point at infinity discussed above.
$\hfill \Box  $

\bigbreak

\noindent It is crucial for the whole study of braids in a handlebody to note that not all
conjugations in the groups $B_{g,n}$  preserve the isotopy class of the closure of an algebraic 
mixed braid.

\begin{cor}{\rm \  Theorem 6 disproves Conjecture 4.4  in \cite{S}.
 } \end{cor}

 \section{Markov functionals}
 
 Theorem 5 opens up the
 possibility to define invariants of links in the handlebody by
 algebraic considerations. This runs largely in parallel with the
 derivation of link invariants in $S^3$ from Markov traces, see for example \cite{J}. In the
 handlebody case, however, trace functionals are not appropriate
 because not all conjugations are allowed, see Theorem 6. Hence, we have to modify
 the definitions, so as to take this into account.
 
\bigbreak

 For an integral domain  $R$ \ let $RB_{g,n}$ denote the group ring
 of the handlebody braid group. A {\it Markov functional} is an 
 $R$-linear map 

\[\mu:RB_{g,n}\longrightarrow R,\] 
 
\noindent for which units $x,\lambda\in R^\ast$ exist such that:
 \begin{eqnarray}
 \mu(1)&=&1\\
 x\mu(\beta)&=&\mu(\iota(\beta))\\
 \mu(\beta\sigma_i^{\pm 1})&=&\mu(\sigma_i^{\pm 1}\beta), \qquad \beta\sigma_i \in B_{g,n} \\
 \mu(\beta_1\sigma_n^{\pm 1}\beta_2)&=&\lambda^{\pm 1}\mu(\beta_1\beta_2) \qquad
 \beta_1,\beta_2\in B_{g,n}
 \end{eqnarray}
 Here $\iota$ stands for the  morphism that embeds $B_{g,n}$ into $B_{g,n+1}$
 by adding an unlinked strand on the right. Moreover, we need the exponent sum of ordinary crossings
 
\[ e: B_{g,n}\longrightarrow \ZZ, \qquad \sigma_i\mapsto 1,
 \ a_r\mapsto 0, \ e(\beta_1\beta_2)=e(\beta_1)+e(\beta_2).\]
 
 \noindent  Theorem 5 now implies the following:  

\begin{defn}{\rm \ The expression defined by

 \[ {\cal L}(L):= x^{n-1}\lambda^{-e(B(L))}\mu(B(L)), \]
 
\noindent where $B(L)\in B_{g,n}$ is an  invariant of oriented links in the handlebody.
  } \end{defn}

 \begin{rem}{\rm \ Nice quotients of the group algebra of $B_{g,n}$ that support
 Markov functionals are to be studied in further work. 
  } \end{rem}

 \bigbreak
 
\noindent {\sc R.H.-O.:  Mathematisches Institut,
G\"{o}ttingen Universit\"{a}t, 

\noindent Bunsenstrasse 3-5, 
D-37073 G\"{o}ttingen, Germany. 

\noindent Homepage: http://www.uni-math.gwdg.de/haering 

\noindent E-mail: roldenburg@gmx.de}

\vspace{.1in}
\noindent {\sc S.L.:  Department of Mathematics,
National Technical University of Athens, Zografou Campus, 
GR-15780 Athens, Greece.  

\noindent Homepage: http://users.ntua.gr/sofial

\noindent E-mail: sofia@math.ntua.gr }

\end{document}